\documentclass[letterpaper]{amsart}
\usepackage{amsfonts}
\usepackage{amssymb}
\usepackage{amsmath}
\usepackage{graphicx}
\usepackage{psfrag}
\usepackage{color}
\usepackage{cite}
\newtheorem{theo}{Theorem}
\newtheorem{coro}[theo]{Corollary}
\newtheorem{prop}[theo]{Proposition}

\newtheorem{lema}[theo]{Lemma}

\newcommand{\R}{\mathbb{R}}

\newcommand{\N}{\mathbb{N}}
\newcommand{\eps}{\varepsilon}

\title{Regularity of  Mediatrices in Surfaces}

\author{Pilar Herreros, Mario Ponce and  J.J.P. Veerman}

\begin{document}

 \begin{abstract}
 For distinct points $p$ and $q$ in a two-dimensional Riemannian manifold, one defines their mediatrix $L_{pq}$ as the set of equidistant points to $p$ and $q$.  It is  known that mediatrices have a cell decomposition consisting of a finite number of branch points connected by Lipschitz curves. This paper establishes additional geometric regularity properties of mediatrices. We show that mediatrices have the radial linearizability property, which implies that at each point they have  a geometrically defined derivative in the branching directions. Also, we study the particular case of mediatrices on spheres, by showing that they are Lipschitz simple closed curves exhibiting at most countably many singularities, with finite total angular deficiency.
 \end{abstract}
 \maketitle

\section{Introduction}

Let $M$ be a compact, connected Riemannian manifold. For any distinct points $p, q\in M$ the {\em mediatrix $L_{pq}$} is defined as the set of points with equal distance  to $p$ and $q$, that is
\[
L_{pq}=\left\{x\in M \ | \ d(x, p)=d(x, q)\right\},
\]
where $d(\cdot, \cdot)$ is the Riemannian distance in $M$.  These geometrical objects can be found in the very beginning of the foundations of geometry. As an illustration, Leibniz proposed to define {\em planes} as mediatrices  in the space. Mediatrices appear in the literature with multiple names, and have been developed by many communities, often separately. The first non-trivial approach to mediatrices appears in the book \cite{BUS} by Busemann, under the name of {\it bisectors}. There, the last chapter is mainly devoted to characterize metric spaces with flat bisectors. The remarkable  result by Busemann is that metric spaces having flat bisectors are (up to a cover) isometric to a finite dimensional euclidian, hyperbolic, or spherical space. In a slightly more general framework, where one replaces points $p, q$ by disjoint compact sets $P, Q\subset M$, the locus of points with equal distance  to $P$ and $Q$ have been studied using the names {\em conflict sets} (see for instance \cite{Yom}, \cite{BS}), {\em midsets} (see \cite{LOV}), {\em equidistant sets} (see \cite{WIL}, \cite{PONSAN}), and {\em skeletal or medial structures} (see \cite{DAM}). \\

The study of mediatrices came about to answer some questions that
came up in the study of {\em Focal Decomposition} pioneered by
\cite{PT}. This point of view was taken in \cite{mediatrix1}.
In spite of the work done in \cite{PT}, there still are many open
problems related to even the simplest cases such as that of the
``mathematical pendulum".
The simple observation of that mediatrices are level sets of the difference of the distance functions to the corresponding points, allows to  consider the  relationship between the cut locus and mediatrices. Namely,  the distance function to a point is differentiable outside the cut locus and hence, the Implicit Function Theorem yields that a mediatrix is differentiable outside the union of the cut locus of the underlying points.  It is well-known (Theorems by Singer-Gluck and Itoh, see for instance \cite{GKM}, \cite{doC}) that the cut locus even for  a $C^{\infty}$ surface of revolution can be a pretty awful set. In contrast the mediatrix on any $C^r, r\geq 3$ surface is locally really surprisingly well-behaved (as we show below).\\

For $p\neq q\in M$, the mediatrix is the boundary between the  two disjoint connected components $\left\{d(x, p)<d(x, q)\right\}$ and $\left\{d(x, q)<d(x, p)\right\}$ of $M\setminus L_{pq}$.
Notice that  by removing any single point of the mediatrix  these sets stop being separated. Hence one says that mediatrices are {\em Minimally Separating Sets}. Separating Sets and especially Minimally Separating Sets are
interesting since they occur in a variety of places throughout
mathematics. In the first place they occur as
 boundaries between the regions where different types of
 Focalization occurs as explained in \cite{PT} and
\cite{mediatrix1}. The classical `Lakes of Wada' construction is
another example (see \cite{HY}). Basin boundaries may in some
cases be another instance of where these sets naturally occur.
 The latter are related to Julia sets, Birkhoff attractor and so on. Recently
other topological properties of this kind of sets have gained some
attention, namely the fact that some are indecomposable continua
(see for example \cite{SKGY}). \\

In \cite{mediatrix1}, Bernhard and Veerman show  that mediatrices in surfaces have a cell decomposition consisting of a finite number of branch points connected by Lipschitz curves. This paper establishes additional geometric regularity properties of mediatrices in surfaces (see definitions below). \\

\noindent{\bf Main Theorem (short version)}. {\em Mediatrices in surfaces have the radial linearizability property.}\\

This  implies that at each point mediatrices  have a geometrically defined derivative in the branching directions.  We devote Section \ref{spheres} to the study of the particular case of mediatrices on spheres, by showing that they are Lipschitz simple closed curves exhibiting at most countably many singularities, with finite total angular deficiency.

\section{Preliminaries}
In what follows, $M$ is a two-dimensional compact Riemannian manifold of class $C^r, r\geq 3$. \\
\paragraph{\bf Minimizing Geodesics. }
Let $x\in M$ and $\theta$ in the unit tangent sphere $S_xM$. We write $\gamma_{x, \theta}:\R\to M$ the unit speed geodesic satisfying $\gamma_{x, \theta}(0)=x$ and $\gamma_{x, \theta}'(0)=\theta$. For $t^*>0$ we say that $\gamma_{x, \theta}$ is a minimizing geodesic between $x$ and $\gamma_{x, \theta}(t^*)$ if the length of
$\gamma_{x, \theta}$ equals the infimum of lengths among absolutely continuous curves joining $x$ to $\gamma_{x, \theta}$, that is, $length(\gamma_{x, \theta}\big |_{[0, t^*]})=d(x, \gamma_{x, \theta}(t^*))$. The triangle inequality yields that $\gamma_{x, \theta}$ is minimizing for every $0<t<t^*$. \\

\paragraph{\bf Polar Coordinates. }
Recall that the exponential map  $exp_x$  at $x\in M$ is defined as
\[
exp_x(v)=\gamma_{x, \hat v}(|v|) \textrm{, where }v\in T_xM \textrm{ , }\hat v=\frac{v}{|v|}\textrm{ , and }exp_x(0)=x.
\]
For small $\rho>0$ (less than the {\em injectivity radius $injrad(M)>0$}), one has that $exp_x\big |_{\{|v|<\rho\}}$ is a diffeomorphism.
By fixing a polar coordinates system in $T_xM$, one can define (via $exp_x$) a local polar coordinates system centered at $x$ on $B_{x, \rho} =\left \{ z\in M\ | \ d(x, z)<\rho\right \}$. In that way,  a point $z\in  B_{x, \rho}, z\neq x $, is identified with a pair $(\theta, d)\in \mathbb{S}^1\times (0, \rho)$, such that $z=\gamma_{x, \theta}(d)$, $d=d(x, z)$ (recall that geodesics are locally minimizing curves, see for instance \cite{doC}).
\\

\paragraph{\bf Minimizing Directions. }
For $p\in M$, $\Theta_{x, p}\subset S_xM$ is the set of directions $\theta$ such that $\gamma_{x, \theta}$ is a minimizing curve from $x$ to $p$. We endow $S_xM$ with the usual arc-distance in the circle. Since $M$ is complete one has
\begin{lema}[Hopf-Rinow Theorem, see \cite{doC}]
For $x\neq p$ the set $\Theta_{x, p}$ is non-empty and compact. $\square$
\end{lema}
\begin{lema}[see \cite{mediatrix3}]
Let $p\neq q\in M$ and $x\in L_{pq}$. Then $\Theta_{x, p}\cap \Theta_{x, q}=\emptyset$.
\end{lema}
\noindent{\bf Proof.} Let $\theta\in \Theta_{x, p}\cap \Theta_{x, q}$ and
$d=d(x,p)=d(x,q)$. Then $\gamma_{x, \theta}(d)=q=p$, which is impossible. $\square$\\
\begin{lema}
Let $x\neq p$,  and two convergent sequences $x_n\to x$, $\theta_n\to \theta$, with $\theta_n\in \Theta_{x_n, p}$ for every $n $ (we identify $S_{x_n}M$ with $S_xM$ via a local chart).  Then one has $\theta\in \Theta_{x, p}$.
\end{lema}
\noindent{\bf Proof.} Let $d_n=d(x_n,p)$. Continuity of the distance function implies $d_n\to d(x,p)$. We have $\gamma_{x_n, \theta_n}(d_n)=p$ for every $n$. Taking the limit we obtain $\gamma_{x, \theta}(d(x,p))=p$, and hence $\theta\in \Theta_{x, p}. \square$ \\

\begin{coro}\label{rhoys}
If $x\in L_{pq}$ then there exist $\rho_x>0$ and $\beta_x\in (0, \pi]$ such that
\[
d_{S_zM}(\Theta_{z, p}, \Theta_{z, q})\geq \beta_x
\]
for all $z\in B_{x, \rho_x}$.$\square$
\end{coro}

\noindent{\bf Proof.}
Given $x$ we know from Lemma 2 that $\Theta_{x, p}\cap \Theta_{x, q}=\emptyset$. Assume by contradiction that there exist sequences $z_n\to x$, $\theta_{p, n}\in \Theta_{z_n, p}$, $\theta_{q, n}\in \Theta_{z_n, q}$ with $|\theta_{p, n}-\theta_{q, n}|\to 0$. Choosing subsequences we can assume that there exists $\theta\in S_{x}M$ such that $\theta_{p, n}\to \theta$, $\theta_{q, n}\to \theta$. Lemma 3
implies that $\theta\in \Theta_{x, p}\cap \Theta_{x, q}$, which is impossible. $\square$\\

For $p\neq q\in M$ and $x\in L_{pq}$, fix a positive orientation in $S_xM$. A closed
connected set $P$ in $S_xM$ is called a {\it pre-wedge at $x$} if its end-points are minimizing directions to $p$ and $q$ respectively, and there are no other minimizing directions to $p$ or $q$ in its interior. Such a set is clearly an interval
in the circle, and with a slight abuse of notation we denote it as:
$P=[\theta_p,\theta_q]\subset S_xM$ (or $[\theta_q,\theta_p]$).
Thus $P$ satisfies:
 \begin{itemize}
\item $\theta_p\in \Theta_{x, p}$, $\theta_q\in \Theta_{x, q}$.
\item $[\theta_p, \theta_q]\cap \left(\Theta_{x, p}\cup \Theta_{x, q}\right)=\left\{\theta_p, \theta_q\right\}$.
\end{itemize}
The \emph{midpoint} $\hat \theta$ of a wedge $P=[\theta_p,\theta_q]$ is the
unique point \emph{contained in the pre-wedge} that equidistant from its endpoints.

\begin{lema}[see \cite{mediatrix3}]
There are finitely many pre-wedges at $x\in L_{pq}$.
\end{lema}
\noindent{\bf Proof. } Assume $\left([\theta_{p}^n, \theta_q^n]\right)_{n\in \N}$ is an infinite sequence of different pre-wedges at $x$. We must have $\lim_{n} |\theta_p^n-\theta_q^n|=0$, since pre-wedges have disjoint interiors. We conclude that the sequences of end-points $(\theta_p^n)_{n\in \N}$, $(\theta_q^n)_{n\in \N}$ have a common accumulation point $\hat \theta$. Compactness of $\Theta_{x, p}$ and $\Theta_{x, q}$ yields that $\hat \theta$ belongs to $\Theta_{x, p}\cap \Theta_{x, q}$ which is impossible. $\square$
\\

For $0<\rho<injrad(M)$ and a pre-wedge $[\theta_p, \theta_q]\subset S_xM$,  we define the corresponding {\it wedge at $x$ of radius $\rho$} as
\[
W_{x, \rho}[\theta_p, \theta_q]=\left\{exp_x(r, \theta)\ | \ 0\leq r<\rho,\ \theta\in [\theta_p, \theta_q] \right\}\subset M.
\]
In \cite{mediatrix3} the authors show how the mediatrix $L_{pq}$ is located in $M$ with respect to the wedges at $x\in L_{pq}$ in the following sense:

\begin{theo}[Bernhard and Veerman \cite{mediatrix3}] Let $M$ be a 2-dimensional compact Riemannian manifold of class $C^r, r\geq 3$. For every pair of distinct points $p, q\in M$ the mediatrix $L_{pq}=\{x\in M \ | \ d(x, p)=d(x, q)\}$ verifies
\begin{enumerate}
\item Let $x\in L_{pq}$. If the cardinalities $\sharp \Theta_{x, p}$ and $\sharp \Theta_{x,q}$ are both 1, then for small $\rho$ the intersection $L_{pq}\cap B_{x, \rho}$ is a continuous simple curve passing through $x$ and differentiable at $x$.
\item Let $x\in L_{pq}$ and $0<\rho<injrad(M)$. The intersection $L_{pq}\cap B_{x, \rho}$ is contained in the finite union of the wedges of radius $\rho $ at $x$.
\item Let $W_{x, \rho}[\theta_p, \theta_q]$ be a wedge at $x\in L_{pq}$. The intersection $L_{pq}\cap W_{x, \rho}[\theta_p, \theta_q]$ is a Lipschitz simple curve connecting $x$ to  a boundary point of  $W_{x, \rho}[\theta_p, \theta_q]$. This curve is called a {\it spoke} at $x$.
\item There are finitely many points $x\in L_{pq}$  so that there exists more than two pre-wedges at $x$.
\item The mediatrix $L_{pq}$ is homeomorphic to a finite and closed simplicial 1-complex.
\end{enumerate}
\end{theo}

Notice that $1.$ above is a direct consequence of the fact that the distance function to a fixed point  $p\in M$ is differentiable outside the cut locus of $p$ and a direct application of the Implicit Function Theorem. The Theorem above says that $L_{pq}\cap  B_{x, \rho}$ is the finite union of two or more Lipschitz spokes emanating from $x\in L_{pq}$. The goal of this article is to study the behavior of a single spoke as it approaches $x\in L_{pq}$. We show that spokes are {\it radially linearizable} in the following sense:\\

\noindent{\bf Definition. } Let $\phi:[0,1]\to M$ be a continuous simple curve. We say that $\eta$ is {\it radially linearizable} at $x=\phi(0)\in M$  if there exists a direction $\hat \theta\in S_xM$ such that
\[
\lim_{t\to 0^+} \theta(t)=\hat \theta,
\]
where $\phi(t)=\left(\theta(t), d(t)\right)$ in polar coordinates centered at $x$.\\

\noindent{\bf Main Theorem [Radial Linearizability of Mediatrices]. }{\em
Let $M$ be a 2-dimensional compact Riemannian manifold of class $C^r, r\geq 3$. Let $p\neq q \in M$,  $x\in L_{pq}$, and $W_{x, \rho}[\theta_p, \theta_q]$ a wedge at $x$. The spoke $L_{pq}\cap W_{x, \rho}[\theta_p, \theta_q]$ is radially linearizable at $x$ in the direction of the midpoint $\hat \theta$ of $[\theta_p, \theta_q]$. }\\

Notice that this Theorem asserts that $L_{pq}\cap W_{x, \rho}[\theta_p, \theta_q]$ is tangent to the bisector of the angle $\angle \theta_p\theta_q$, and hence, $L_{pq}$ is tangent to the directions of the equidistant set to $\Theta_{x, p}$ and $\Theta_{x, q}$ in $S_xM$. A similar result was conjectured in \cite{mediatrix3}  (compare also with \cite{BS}).\\

\noindent{\bf Corollary. } {\em Mediatrices in surfaces have no cusp like points.}

\section{Proof of the Main Theorem}
\paragraph{\bf The distance function has one-side directional differential. } For $p\in M$,
 the triangle inequality implies that the function $x\mapsto d(x,p)$ is $1$-Lipschitz. Even though $d(x,p)$ is not necessarily differentiable in $M$ (in fact, it is not differentiable at the cut locus of $p$, see for instance \cite{doC}), for every $x\in M$ this function has a one-sided derivative in every direction. More precisely, let $\theta \in S_xM$. The limit
\[
\lim_{t\to 0^+}\frac{d(\gamma_{x, \theta}(t),p)-d(\gamma_{x, \theta}(0),p)}{t}
\]
does exist. This fact, and the explicit formula for the limit, is apparently a well known fact that can be deduced from the first variation formula. Nevertheless, just recently, and due to the absence of any unambiguous literature in the subject,  I. Adelstein \cite{Adel} has produced an explicit proof (see also Section 4.5 in \cite{BBI}).
\begin{prop}\label{folk}
Using the above notation, one has
\[
\lim_{t\to 0^+}\frac{d(\gamma_{x, \theta}(t),p)-d(\gamma_{x, \theta}(0),p)}{t}
= -\cos(d_{S_xM}(\theta, \Theta_{x,p}))
\]
where $d_{S_xM}(.,.)$ is the Riemannian distance on the unit circle
between $\theta$ and the compact set of minimizing directions from $x$ to $p$.
$\square$

\end{prop}
\begin{lema}\label{l.1234}
Let $g:[0, 1)\to \R$ be a continuous function and  $c\in R$. Assume that the one-side derivative exists and verifies
\[
g'_+(t) \equiv \lim_{s\to 0^+}\frac{g(t+s)-g(t)}{s}>c
\]
for every $t\in [0, 1)$. Then $g(t)>g(0)+ct$ for every $t\in (0, 1)$.
\end{lema}
\noindent{\bf Proof. }  By considering the function $g(t)-ct$ we can reduce us to the case $c=0$. For every $t\in [0, 1)$ there exists $\delta_t>0$ such that $g(t+s)-g(t)>0$ for $0<s<\delta_t$. Hence $g$ can not have local maximal points. As $g$ is continuous, the only possibility is that $g$ is strictly increasing. $\square$\\

In what follows we use the notation of the statement of the Main Theorem, that is,  we fix a pre-wedge $[\theta_p, \theta_q]$ at $x\in L_{pq}$,  write $\hat \theta$  the midpoint of $[\theta_p, \theta_q]$, and   $\rho_x>0$ given by Corollary \ref{rhoys}.

\begin{lema}\label{l.1}
If $L_{pq}\cap W_{x, \rho}[\theta_p, \theta_q]\cap \gamma_{x, \theta}([0, \rho_x))$ contains a sequence of points converging to $x$ then $\theta=\hat \theta$.
\end{lema}
\noindent{\bf Proof.} Define $f_{pq}(z)=d(z,q)-d(z,p)$. Hence $L_{pq}=f_{pq}^{-1}(0)$.  Proposition \ref{folk} allows to deduce that $f_{pq}$ has one-side differential at $x$ in the direction of $\theta$  equals to $-cos|\theta-\tilde \theta_p|+\cos|\theta-\tilde \theta_q|$ for some $\tilde \theta_p\in \Theta_{x,p}$ , $\tilde \theta_q\in \Theta_{x,q}$. If the restriction of $f_{pq}$ to the ray $\gamma_{x, \theta}$ has roots accumulating at $x$ then this one-side differential must be equal to zero. The only possible case is $\tilde \theta_p=\theta_p, \tilde \theta_q=\theta_q$ (otherwise minimizing geodesics
will cross) and it follows that the only solution is $\theta=\hat \theta$. $\square$\\

This Lemma says that whenever a spoke crosses infinitely many times a geodesic ray emanating from $x$, and the crossing points approach $x$, then the geodesic ray necessarily points in the direction of the bisector of the corresponding pre-wedge.  The next result estimates for how much the points in the geodesic ray $\gamma_{x, \hat \theta}$  fail to belong to the mediatrix $L_{pq}$.

\begin{lema}\label{l.2}
Given $p$ and $q$, then for each $x\in L_{pq}$ there exists a positive function
$\eta$ such that
\[
\left| d\left(\gamma_{x, \hat \theta}(t),q\right)-
d\left(\gamma_{x, \hat \theta}(t),p\right)  \right |\leq
\eta(t)t\textrm{ for } 0\leq t<\rho ,
\]
and $\lim_{t\to 0} \eta(t)=0$.
\end{lema}

\noindent{\bf Proof.} The function $t\mapsto f_{pq}\left (\gamma_{x, \hat \theta}(t)\right)$ together with its first (right directional) derivative vanish at $t=0$ and the result follows. $\square$ \\

The next result tells us that close to points that fail to belong to $L_{pq}$ by few, we can find a genuine point in the mediatrix. This result is inspired by  Section 4.2 of \cite{PONSAN}.

\begin{lema}\label{l.3}
Let $x\in L_{pq}$ and $\rho_x>0, \beta_x\in (0, \pi)$ given by Corollary \ref{rhoys}. There exists $ \eps_x>0$  such that for all  $z\in B_{x, \rho_x/2}$ that verifies $
\left| d(z,p)-d(z,p)\right|< \eps_x
 $  there exists $z^*\in L_{pq}$, and
\[
d(z, z^*)<\frac{\left |  d(z,q)-d(z,p)\right|}{1-cos(\beta_x)}.
\]
\end{lema}

\begin{figure}
 \psfrag{r}{$\rho_x$}
 \psfrag{p}{$p$}
  \psfrag{q}{$q$}
   \psfrag{z}{$z$}
  \psfrag{z*}{$z^*$}
   \psfrag{x}{$x$}
  \psfrag{L}{$ L_{pq}$}
     \includegraphics[scale=0.8]{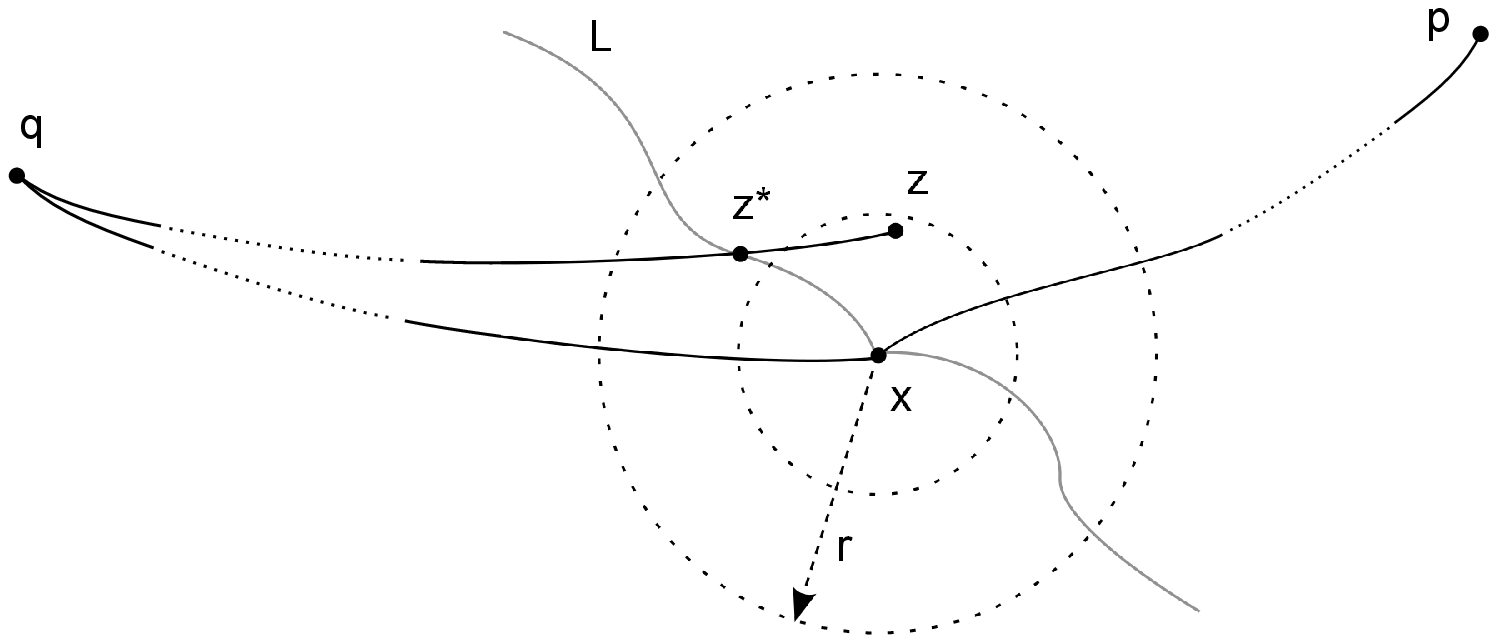}
\end{figure}

 \noindent{\bf Proof. } Let's assume that $d(z,q)>d(z,p)$. Pick $\theta\in \Theta_{z, q}$ (this also gives $\theta\notin \Theta_{z, p}$). By moving continuously along the geodesic ray $z(t)\equiv\gamma_{z, \theta}(t)$ we must encounter a point $z^*\equiv z(t^*)$ in the mediatrix $L_{pq}$ (since $\gamma_{z, \theta}$ reaches $q$). We estimate the distance $t^*$ between $z$ and $z^*$. Note that since $d(x,z)<\rho_x/2$ then by construction $d(x,z(t))<\rho_x$ at least for $t\leq \rho_x/2$.

By Proposition \ref{folk} we have
\[
\lim_{t\to 0^+}\frac{d(z(t),p)-d(z(0),p)}{t}
= -\cos(d_{S_{z(t)}M}(z'(t), \Theta_{z(t),p}))
\]
For every $t\in [0, \rho_x/2]$ one has $z'(t)\in \Theta_{z(t), q}$ and so Corollary
\ref{rhoys} implies that there is a $\beta_x\in(0,\pi)$ so that
\[
 \cos(d_{S_{z(t)}M}(z'(t), \Theta_{z(t),p}))\leq \cos \beta_x
\]
Now Lemma \ref{l.1234} implies that for every $t\in (0, \rho_x/2]$
 \[
 d(z(t),p)\geq d(z,p)-t \cos\beta_x.
 \]
 We also have that
 \[d(z(t),q)=d(z,q)-t.
 \]
 Putting these two together we have that
 \[
 \begin{array}{ccl}
 d(z(0),q)-d(z(0),p)&>& 0 \\
 d(z(t),q)-d(z(t),p) & \leq & d(z(0),q)-d(z(0),p)-t(1-\cos(\beta_x))
 \end{array}
 \]
The latter is less than or equal to zero at $t_0=\frac{d(z,q)-d(z,p)}{1-\cos(\beta_x)}$. By continuity of the distance function, we see that
 \[
 t^{*}\leq\frac{d(z,q)-d(z,p)}{1-\cos(\beta_x)}.
 \]

If we choose $\eps_x \equiv \frac{1}{2}\rho_x(1-\cos(\beta_x))>0$, so that $t_0<\rho_x/2$, we stay inside $B_{x, \rho_x}$  and the Lemma is verified. $\square$
 \\

\noindent{\bf Proof of the Main Theorem. }Pick $x\in L_{pq}$ and a wedge  $W_{x, \rho}[\theta_p, \theta_q]$.  Reasoning by contradiction let's suppose that there exists $\tau>0$ such that the spoke contains a sequence of points $x_n$ converging to $x$ outside the cone $[\hat \theta -\tau, \hat \theta +\tau]\times (0, \rho)$ (in polar coordinates centered at $x$). Lemmas \ref{l.2} and \ref{l.3} assert that the spoke contains a sequence of points inside the cone as it approaches $x$, since the spoke passes trough a point at a distance less than  $ (1-\cos\beta_x)^{-1}{\eta(t)t}$ from  $\gamma_{x, \hat \theta}(t)$ and thus inside the cone for $t$ sufficiently small. By the continuity of the spoke, it must cross $ \gamma_{x, \hat \theta+\tau} $  (or $ \gamma_{x, \hat \theta-\tau} $) an infinite number of times, and Lemma \ref{l.1} says that this is impossible unless $\tau=0$. $\blacksquare$

\section{Mediatrices in spheres, the Pugh-Tangerman example}\label{spheres}
Let $M$ be a simply connected 2-dimensional Riemannian manifold of class $C^r, r\geq 3$ (a {\em sphere}). For $p\neq q \in M$ we already know that there are at most finitely many points $x\in L_{pq}$ such that $L_{pq}$ is not locally homeomorphic to a real interval at $x$ (that is, there are more than two wedges at $x$). Let us call such a point a {\em branching point} (and a {\em simple point} to a non-branching point). For a simple point $x$ the Main Theorem asserts that $L_{pq}$ has two well defined directions $\hat \theta_1, \hat \theta_2\in S_xM$ at $x$ (the bisectors of the two pre-wedges at $x$, or in other words, the
equidistant set in $S_xM$ of the sets of minimizing directions from $x$ to $p$ and $q$ respectively). We define the {\it deficiency angle at $x$} as
\[
\textrm{def}_{L_{pq}}(x)=\left |\pi - |\hat \theta_1-\hat \theta_2| \right |.
\]
If $\textrm{def}_{L_{pq}}(x)\neq 0$ we say that $x$ is a {\em Lipschitz singularity} of the mediatrix $L_{pq}$.\\

\begin{figure}
 \psfrag{p}{$p$}
  \psfrag{q}{$q$}
     \includegraphics[scale=0.4]{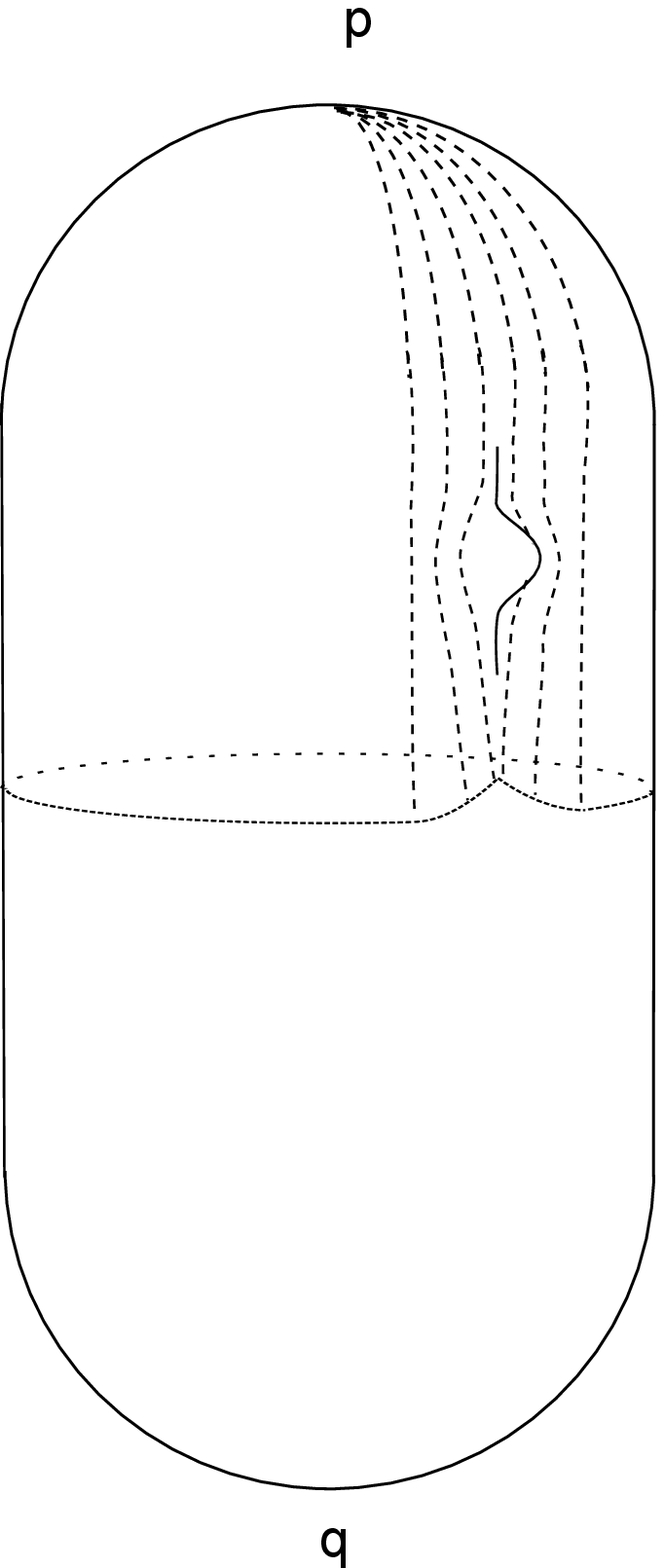}
\end{figure}

\noindent{\bf Pugh-Tangerman example.} A very illustrative example is what has
come to be known as \emph{Pugh's Cigar $C$} embedded in $\R^3$.
It consists of the cylinder given by  $x^2+y^2=1$ plus two (unit) hemispheres
glued to each of the boundaries. Let $p$ be its North Pole $(0,0,1)$ and $q$
its South Pole. The mediatrix is given by the equator at $z=0$.

We can obtain an arbitrarily smooth cigar by smoothing around $|z|=\pm 1$.
If we add a ``bump" in the flat ``northern" part of the cigar, keeping $C$ as a smooth manifold, is easy to see that the mediatrix will develop a Lipschitz singularity near the closest point to the bump. Adding a sequence of bumps
in the northern part that do not intersect each other we can create a cigar whose mediatrix will have countably many Lipschitz singularities.

In what follows we want to show that this is, in some sense, the worst possible
behavior for mediatrices on spheres. This proposition is inspired on conversations
in the 1990's by Pugh and Tangerman with Veerman at Rockefeller
University.

\begin{prop}
Let $M$ be a sphere of class $C^r, r\geq 3$ and $p\neq q \in M$. Then the
mediatrix $L_{pq}$ verifies
\begin{enumerate}
\item $L_{pq}$ is a Lipschitz simple closed curve (there are no branching points).
\item There are  at most countably many Lipschitz singularities in $L_{pq}$. Moreover the deficiency angles form an absolutely summable sequence, that is
\[
\sum_{x\in L_{pq}}\textrm{def}_{L_{pq}}(x)<\infty.
\]
\end{enumerate}
\end{prop}
\noindent{\bf Proof. } Simple closed follows from earlier results
(\cite{mediatrix3}). The novelty here is the Lipchitz property.
The fact that the mediatrix is a simple closed curve implies that at each
$x\in L_{pq}$ there are exactly two spokes, and thus exactly two wedges.
By our Main Theorem, the mediatrix is tangent to the bi-sectors of these
two wedges, and  thus is Lipschitz. This proves (1).

We saw that at a given point $x\in L_{pq}$ there are exactly two wedges
each of which is bounded by two minimizing geodesics $\gamma_{x,p}$ and
$\gamma_{x,q}$. Take the complement of the closure of these wedges.
If it is the empty set, the mediatrix is differentiable at $x$. If not,
then (a) that complement is a maximal wedge shaped set bounded by two
minimizing geodesics from $x$ to $p$, or (b) bounded two minimizing geodesics
geodesics from $x$ to $q$, or (c) both of these occur.
In case (a) let us denote by $\mathcal{J}_{x, p}$ the open domain
bounded by the indicated geodesic rays from $x$ to $p$ and by
$\mu_{x,p}$ the opening angle at $x$ of that wedge-shaped figure. Similarly we have
$\mathcal{J}_{x, q}$ and $\mu_{x,q}$ for case (b).
Case (c) means that both $\mathcal{J}_{x, p}$ and $\mathcal{J}_{x, p}$
are non-empty. Notice that $\mathcal{J}_{x, p}=\emptyset$ if and only if $\sharp \Theta_{x, p}=1$.

If we have two minimizing geodesics form $x_1, x_2\in L_{pq}$ to $p$ or $q$, by triangle inequality the intersection point would have to be strictly closer to the first endpoint than the second, thus by symmetry these geodesics do not intersect. Therefore the collection (as $x$ runs on $L_{pq}$) of Jordan domains $\mathcal{J} _{x, p}, \mathcal{J}_{x, q}$ are pairwise disjoint. We conclude that there are at most countably many points $x \in L_{pq}$ such that $\sharp \Theta_{x, p}$ or $\sharp \Theta_{x, q}$ are different from $1$. A direct examination yields that
\[
2\textrm{def}_{L_{pq}}(x)=\big| \mu_{x, p}-\mu_{x, q}\big|
\]
The Gauss-Bonnet Theorem applied to $\mathcal{J}_{x, p}$ gives
\[
\int_{\mathcal{J}_{x, p}}k(w)dw=\mu_{x,p}+\alpha_{x, p},
\]
where $\alpha_{x, p}$ is the angle of the landing geodesic rays make at $p$
(note that $mu_{x,p}$ is the angle they make at $x$),
and $k(w)$ is the curvature at $w\in M$. We estimate
\begin{eqnarray*}
\sum_{x\in L_{pq}}\textrm{def}_{L_{pq}}(x)&<& \frac{1}{2}\sum_{x\in L_{pq}}
\left(\int_{\mathcal{J}_{x, p}}|k(w)|dw+\int_{\mathcal{J}_{x, q}}|k(w)|dw+|
\alpha_{x, p}|+|\alpha_{x,q}|\right)\\
&\leq& \frac{1}{2}\left(\int_{M}|k(w)|dw+2\pi\right).\square
\end{eqnarray*}
\begin{small}

\vspace{0.7cm}
\noindent{\bf Acknowledgments.} Pilar Herreros acknowledges the financial support
of the research project FONDECYT 11121125. Mario Ponce and J. J. P. Veerman
acknowledge the financial support of the research project FONDECYT 1140988.
Mario Ponce received partial financial support from the project PIA-CONICYT
ACT-1103. J.J.P. Veerman's research was also partially supported by the European
Union's Seventh Framework Program (FP7-BEGPOT-2012-2013-1) under grant agreement
n316165.

\end{small}


\begin{footnotesize}

\vspace{0.25cm}

\noindent{Pilar Herreros}

\noindent{Facultad de Matem\'aticas, Pontificia Universidad Cat\'olica de Chile}

\noindent{e-mail: pherrero@mat.puc.cl}\\

\noindent{Mario Ponce}

\noindent{Facultad de Matem\'aticas, Pontificia Universidad Cat\'olica de Chile}

\noindent{e-mail: mponcea@mat.puc.cl}\\

\noindent{J.J.P. Veerman}

\noindent{Department of Mathematics and Statistics, Portland State
University, and} \\
\noindent{CCQCN, Dept. of Physics, University of Crete, 71003 Heraklion, Greece}

\noindent{e-mail: veerman@pdx.edu}\\

\end{footnotesize}
\end{document}